\documentclass[a4paper,english,brazil]{article}
\usepackage[T1]{fontenc}
\usepackage[latin1]{inputenc}
\usepackage{amsmath}
\usepackage{amssymb}

\makeatletter

\@ifundefined{date}{}{\date{}}

\usepackage{amsthm}\usepackage{amsfonts}\usepackage{amstext}

\makeatother

\usepackage{babel}
\begin{document}
\selectlanguage{english}%

\title{A New Algorithm for Evaluating Determinants \\ -- The Reduction
Method}

\author{Ricardo S. Vieira%
\thanks{E-mail: rsvieira@df.ufscar.br.%
} }

\maketitle
\begin{center} \textit{Universidade Federal de São Carlos, Departamento
de Física, \\ Caixa Postal 676, CEP 13569-905, São Carlos, Brazil
}\end{center}

\section*{Abstract}
\begin{quotation}
We present here a new method for evaluating determinants -- the reduction
method%
\footnote{The relationship between this work and that of the mathematicians
C. L. Dodgson and F. Chiò is presented in the Addendum.%
}. Firstly, in the section 2, we apply it to third-order determinants
and after, in the section 3, we generalize it to higher-order determinants.
In the section 4 an alternative formulation of the reduction method
is presented and, in the section 5, we give the respective mathematical
proofs. 
\end{quotation}

\section{Introduction}

It is not necessary to say that theory of determinants is very important
both in physics as in mathematics. Effectively, in several circumstances,
the solution of a given problem falls on the resolution of a determinant.
If we need to solve a second or third-order determinant, we have no
problem at all, since the calculations are pretty easy in these cases.
Fourth-order determinants can also be easily solved by the well-known
Laplace method, so we are carefree in this case too. However, if the
determinant is of fifth or higher-order, then we begin to meet with
hand-computational problems, since the resolution of such determinants,
even with the Laplace method, could be a hard and boring job. The
use of a computer machine is generally indispensable in these cases.

To attack this problem, we shall present here a new algorithm for
determinant evaluation, which we believe to be the simplest method
for hand-evaluation of higher-order determinants. Effectively, to
use this method the reader needs only to know how evaluate determinants
of order 2. At each step of its application, the order of the determinant
it is reduced by one, which justifies the name of reduction.

We hope that the simplicity of this method may contribute to a better
presentation of the theory of determinants, even in elementary levels.
In fact, the exposition of the theory at this level usually stops
on the theory of fourth-order determinants, due to the technical difficulties
commented above. The reduction method might, thus, eliminate such
difficulties.

The paper is organized as follows. In the section 2 we present the
reduction method for third-order determinants. In the next section
we extend it to higher-order determinants. In the section 4 an alternative
formulation of the reduction method it is presented. Finally, the
mathematical proofs are presented in the last section.

\section{The Reduction Method for Third-Order Determinants}

Let $A$ be a general third-order determinant: 
\begin{equation}
A=\left|\begin{array}{ccc}
a_{11} & a_{12} & a_{13}\\
a_{21} & a_{22} & a_{23}\\
a_{31} & a_{32} & a_{33}
\end{array}\right|.
\end{equation}

The first step of the reduction method consists in evaluate the four
minors $b_{ij}$, built up from the adjacent elements of $A$, namely
\begin{equation}
b_{11}=\left|\begin{array}{cc}
a_{11} & a_{12}\\
a_{21} & a_{22}
\end{array}\right|,\hspace{0.75cm}b_{12}=\left|\begin{array}{cc}
a_{12} & a_{13}\\
a_{22} & a_{23}
\end{array}\right|,\hspace{0.75cm}b_{21}=\left|\begin{array}{cc}
a_{21} & a_{22}\\
a_{31} & a_{32}
\end{array}\right|,\hspace{0.75cm}b_{22}=\left|\begin{array}{cc}
a_{22} & a_{23}\\
a_{32} & a_{33}
\end{array}\right|.
\end{equation}
 This allows us to build the reduced determinant $B$, which is constructed
by disposing the minors $b_{ij}$ in a natural order: 
\begin{equation}
B=\left|\begin{array}{cc}
b_{11} & b_{12}\\
b_{21} & b_{22}
\end{array}\right|.
\end{equation}
 Then, we argue that the evaluation of the determinant $B$ divided
by $a_{22}$ (i.e., divided by the central element of $A$), results
to be the same as the determinant $A$, provided that the element
$a_{22}$ be different from zero: 
\begin{equation}
A=B/a_{22}.
\end{equation}

The proof is nothing but a matter of calculation. It can be done substituting
the expression of the minors $b_{ij}$ in the expression of $B$ and
doing the calculi explicitly. We left this job to the reader.

Notice, however, that if the central element of $A$ is null, the
reduction method cannot be directly applied. Effectively, in this
case we would meet with a division by zero. Fortunately, this is not
a problem at all. We may, for instance, before applying the method,
permutate cyclically the lines (or columns) of $A$ in such a way
that the new central element turns to be a non-null element. These
operations, as we know, do not change the value of a determinant.

\section{The Reduction Method for Higher-Order Determinants}

Let us now show how we can generalize the reduction method to higher-order
determinants. Let $A=|a_{ij}|$, $1\leq i,j\leq n$, be a $n$-order
determinant. As before, we proceed by building up, from each group
of four adjacent elements of $A$, a total of $n-1$ minors of second-order,
$b_{ij}$. That is, we construct the minors $b_{ij}$ through the
expression 
\begin{equation}
b_{ij}=\left|\begin{array}{cc}
a_{ij} & a_{i,j+1}\\
a_{i+1,j} & a_{i+1,j+1}
\end{array}\right|,\hspace{1cm}1\leq i,j\leq n-1.
\end{equation}
 With the values of these minors disposed in an adequate way, we build
up the reduced determinant $B=|b_{ij}|$, whose order is $n-1$. Then,
we proceed with the reduction of the determinant $B$ itself. Notice,
however, that the reduced minors $c_{ij}$ must be defined now by
the expression 
\begin{equation}
c_{ij}=\frac{\left|\begin{array}{cc}
b_{ij} & b_{i,j+1}\\
b_{i+1,j} & b_{i+1,j+1}
\end{array}\right|}{a_{i+1,j+1}},\hspace{1cm}1\leq i,j\leq n-2.
\end{equation}
 The division by $a_{i+1,j+1}$ can easily be understood: we have
seen in the precedent section that after the evaluation of the reduced
determinant $B$ was necessary to perform a division (by the central
element of $A$). An analogous thing happens here. After the second
reduction, each minor $c_{ij}$ it is the result of the reduction
of a third-order determinant. Consequently, each one of these reduced
minors must be divided by the respective central element of the third-order
determinant whose it was derived.

With the elements $c_{ij}$ evaluated by the expression above, we
can build the second determinant reduced $C=|c_{ij}|$, now of $(n-2)$-order.
We can carry on the reduction method as we did in the last step: we
build up the $n-3$ minors $d_{ij}$ from the adjacent elements of
$C$, remembering that we also must divide each one of them by the
element $b_{i+1,j+1}$. This procedure must be repeated until the
determinant $A$ becomes a number -- the value of the determinant
$A$. Thus, if $S=|s_{ij}|$ is the determinant obtained from $A$
when we have had applied the reduction method $k$ times ($k\geq2$),
and $R=|r_{ij}|$ is the determinant whose reduction gives place to
$S$, then, the elements $t_{ij}$, which comes from the reduction
of $S$, must be given by 
\begin{equation}
t_{ij}=\frac{\left|\begin{array}{cc}
s_{ij} & s_{i,j+1}\\
s_{i+1,j} & s_{i+1,j+1}
\end{array}\right|}{r_{i+1,j+1}},\hspace{1cm}1\leq i,j\leq n-k-2.
\end{equation}

To illustrate the simplicity of the reduction method, we shall give
at this point an example of a fifth-order determinant. The reader
can notice that the calculations are pretty easy and can really be
made mentally. 
\begin{equation}
\left|\begin{array}{rrrrr}
1 & 0 & -2 & 3 & 2\\
-1 & -3 & 2 & -2 & 0\\
-3 & -2 & 2 & -1 & 1\\
-2 & 3 & -1 & 2 & 0\\
0 & -3 & 1 & -1 & -3
\end{array}\right|\Rightarrow\left|\begin{array}{rrrr}
-3 & -6 & -2 & 4\\
-7 & -2 & 2 & -2\\
-13 & -4 & 3 & -2\\
6 & 0 & -1 & -6
\end{array}\right|\Rightarrow\left|\begin{array}{rrr}
12 & -8 & 2\\
-1 & 1 & -2\\
8 & -4 & -10
\end{array}\right|\Rightarrow\left|\begin{array}{rr}
-2 & 7\\
1 & -6
\end{array}\right|=5.
\end{equation}

It is necessarily comment, however, that the reduction method cannot
be directly applied in some cases because it would lead to a division
by zero. This will happen whenever some internal minor of $A$ is
singular%
\footnote{By singular we want mean that its determinant vanishes. Notice also
that a simple element as $a_{ij}$ should be considered as a minor
of order 1.%
}. In some cases, this problem can be repaired by a simple permutation
of the lines (or columns) of the original determinant, as we did in
the precedent section, but generally this is not enough to fix the
problem. To overcome this difficult at all, we will proceed in another
way. Suppose that the determinant $A$ contains in a singular minor
$M$, whose order is $m<n$. We should notice that, if we substitute
one (or more, if necessary) element of $M$, says $m_{ij}$, by a
parameter $\varepsilon$, then $M$ no more will be a singular minor
(because now $\varepsilon$ is a free parameter). This solves the
problem because we eliminate the singularity with this substitution.
So, the reduction method is now quite applicable. Reducing the determinant
$A$ at the end we will have found a expression that depends on $\varepsilon$;
the value of $A$ is simply obtained by restoring the original value
$\varepsilon=m_{ij}$ (sometimes it is easier restore the value of
$\varepsilon$ soon after the division be done).

In the sequence we shall present an example of this. Let $A$ be the
following fifth-order determinant: 
\begin{equation}
A=\left|\begin{array}{rrrrr}
-1 & 0 & -1 & 0 & -2\\
2 & 1 & -2 & -1 & 0\\
-1 & 2 & 1 & -2 & 1\\
1 & 3 & 1 & -2 & -1\\
-1 & 1 & -2 & -2 & 0
\end{array}\right|.
\end{equation}
 The determinant $A$ has the value $A=15$. Notice that the reduction
method cannot be directly applied here because the minor $b_{33}=(a_{33}a_{44}-a_{34}a_{43})$
is singular (in fact, if we try to apply the method, then we would
to meet us with a division by zero in the third reduction). To remove
this singularity we may substitute, for instance, the element $a_{33}$
by a parameter $\varepsilon$, which might a priori assume any value.
After this we can apply the reduction method without problems. Doing
so, we will find, 
\begin{eqnarray}
 & \left|\begin{array}{rrrrr}
-1 & 0 & -1 & 0 & -2\\
2 & 1 & -2 & -1 & 0\\
-1 & 2 & \varepsilon & -2 & 1\\
1 & 3 & 1 & -2 & -1\\
-1 & 1 & -2 & -2 & 0
\end{array}\right|\Rightarrow\left|\begin{array}{rccr}
-1 & 1 & 1 & -2\\
5 & \varepsilon+4 & \varepsilon+4 & -1\\
-5 & 2-3\varepsilon & 2-2\varepsilon & 4\\
4 & -7 & -6 & -2
\end{array}\right|\Rightarrow\left|\begin{array}{rcr}
-9-\varepsilon & 0 & -7-2\varepsilon\\
15-5\varepsilon & \varepsilon+4 & -9-\varepsilon\\
9+4\varepsilon & 2+4\varepsilon & -10-2\varepsilon
\end{array}\right|\Rightarrow\nonumber \\
\nonumber \\
 & \Rightarrow\left|\begin{array}{cc}
-9-\varepsilon & 7+2\varepsilon\\
8\varepsilon-3 & -11-\varepsilon
\end{array}\right|\Rightarrow30-15\varepsilon.
\end{eqnarray}
 Finally, restoring the value $\varepsilon=1$, we get $A=15$, which
is the correct value of $A$. Notice that the substitution $\varepsilon=1$
might be made soon after the third reduction as well. In this example
the difference is insignificant, but to higher-order determinants
it may be most appreciable. The reader could see from this example
that the introduction of the parameter $\varepsilon$ makes the calculations
more boring in general, but at least it is a solution to the problem.
Anyway, in the next section we shall present an alternative form of
the reduction method which solves the difficulty of the divisions
by zero once and for all.

\section{An Alternative Formulation of the Reduction Method}

In the precedent sections we have seen how the reduction method of
determinants works. We saw that the reduction is employed by the evaluation
of the second-order minors that are formed from the adjacent elements
of the original determinant. This choice of the adjacent elements
for the building it is, however, quite arbitrary. In fact, we could
build these minors from non-adjacent elements instead, although the
correct arrangements would be a hard job. That this is possible can
be evidenced by the following argument: suppose we have changed the
lines and columns of the original determinant $A$ to obtain a new
determinant $A'$, which we consider to have the same value of $A$.
Then, if we build the minors from the adjacent elements of $A$, although
they are no longer adjacent in the determinant $A'$, it is evident
that the reduction method will also work. This shows that we can build
minors from non-adjacent elements, although we should be careful with
the configurations which are allowed.

Now we shall show one of these possible arrangements which came to
be important, because of its simplicity and also because it eliminates
at once the problem of the divisions by zero that we met on the precedent
section. This alternative formulation of the reduction method is based
on a different way of building the second-orders minors $b_{ij}$
from the determinant $A$. Instead building them from the adjacent
elements of $A$, we proceed as follows. We fix a specific element
of $A$, says the element $a_{rs}$, and we build all minors only
from this element. That is, we fix a element $a_{rs}$ of $A$ and,
for each element $a_{ij}$ with $i\neq r$ and $j\neq s$, we build
the minors that contain the elements $a_{rj}$ and $a_{is}$, so as
the element $a_{ij}$ stays in the diagonal with $a_{rs}$ and the
other two remain at the sides from them. Let us explain this in a
better way. Consider the following $n$-order determinant 
\begin{equation}
A=\left|\begin{array}{ccccc}
a_{11} & \hdots & a_{1s} & \hdots & a_{1n}\\
\vdots & \ddots & \vdots & \ddots & \vdots\\
a_{r1} & \hdots & a_{rs} & \hdots & a_{rn}\\
\vdots & \ddots & \vdots & \ddots & \vdots\\
a_{n1} & \hdots & a_{ns} & \hdots & a_{nn}
\end{array}\right|.
\end{equation}
 If we fix the non-null element $a_{rs}$ of $A$, the second-order
minors that we must build will be, respectively, 
\begin{eqnarray}
 & b_{11}=\left|\begin{array}{cc}
a_{11} & a_{1s}\\
a_{r1} & a_{rs}
\end{array}\right|,\hspace{0.75cm}\ldots\hspace{0.75cm}b_{1,n-1}=\left|\begin{array}{cc}
a_{1s} & a_{1n}\\
a_{rs} & a_{rn}
\end{array}\right|,\nonumber \\
\nonumber \\
 & \vdots\hspace{2cm}\ddots\hspace{2cm}\vdots\nonumber \\
\nonumber \\
 & b_{n-1,1}=\left|\begin{array}{cc}
a_{r1} & a_{rs}\\
a_{n1} & a_{ns}
\end{array}\right|,\hspace{0.75cm}\ldots\hspace{0.75cm}b_{n-1,n-1}=\left|\begin{array}{cc}
a_{rs} & a_{rn}\\
a_{ns} & a_{nn}
\end{array}\right|.
\end{eqnarray}
 In a general way, these minors can be written in a more compact form
as 
\begin{equation}
b_{ij}=\sigma_{ij}(a_{ij}a_{rs}-a_{is}a_{rj}),\hspace{1cm}1\leq i,j\leq n-1,
\end{equation}
 where $\sigma_{ij}=\frac{(r-i)(s-j)}{|(r-i)(s-j)|}$ it is just a
function that provides the correct sign for each minor $b_{ij}$.

With those minors thus defined, we can reduce the determinant $A$
to obtain the reduced determinant $B$, of $(n-1)$-order: 
\begin{equation}
B=\left|\begin{array}{ccc}
b_{11} & \hdots & b_{1,n-1}\\
\vdots & \ddots & \vdots\\
b_{n-1,1} & \hdots & b_{n-1,n-1}
\end{array}\right|.\label{B}
\end{equation}

We can continue with the reduction process by reducing at its turn
the determinant $B$. For this, we choose a non-null element $b_{rs}$
of $B$ and we build the minors $c_{ij}$ as before, but now we must
divide each minors thus constructed by the element $a_{rs}$ which
we have fixed in the reduction of the determinant $A$. But since
all minors should be divided by the same element, $a_{rs}$, and since
the determinant $C$ should have $n-2$ lines and columns, these divisions
are the same as to divide the whole determinant by $(a_{rs})^{n-2}$.
That is, after the second reduction, we should define 
\begin{equation}
C=\frac{1}{(a_{rs})^{n-2}}\left|\begin{array}{ccc}
c_{11} & \hdots & c_{1,n-2}\\
\vdots & \ddots & \vdots\\
c_{n-2,1} & \hdots & c_{n-2,n-2}
\end{array}\right|.
\end{equation}
 Then we can continue by reducing the determinant $C$ itself: we
fix a non-null element $c_{rs}$, build the $n-3$ minors $d_{ij}=\sigma_{ij}(c_{ij}c_{rs}-c_{is}c_{rj})$
and divide the resultant determinant by $(b_{rs})^{n-3}$, obtaining
thus the reduced determinant $D$. This process must be repeated until
the value of the determinant $A$ be evaluated. In short, we found
that the determinant $A$ will be related to the last reduced determinant
$Z$ (of order 1), by the expression 
\begin{equation}
A=\frac{Z}{(a_{rs})^{n-2}(b_{rs})^{n-3}\ldots(x_{rs})^{1}},
\end{equation}
 where $x_{rs}$ it is the element fixed at the reduced determinant
$X$, of order 3. We can also notice that, since we should have for
the determinant $B$ an analogous expression, namely, 
\begin{equation}
B=\frac{Z}{(b_{rs})^{n-3}(c_{rs})^{n-4}\ldots(x_{rs})^{1}},
\end{equation}
 follows that, 
\begin{equation}
A=\frac{B}{(a_{rs})^{n-2}}.
\end{equation}
 Hence, the reduction method presented here can be considered as a
recursive application of this last formula.

To illustrate this alternative way of reduction, let us to apply it
to the first example of the precedent section (the fixed elements
in every determinant were written in bold type and the divisions were
left to the final of the reduction process): 
\begin{equation}
\left|\begin{array}{rrrrr}
1 & 0 & -2 & 3 & 2\\
\mathbf{-1} & -3 & 2 & -2 & 0\\
-3 & -2 & 2 & -1 & 1\\
-2 & 3 & -1 & 2 & 0\\
0 & -3 & 1 & -1 & -3
\end{array}\right|\Rightarrow\left|\begin{array}{rrrr}
-3 & 0 & 1 & 2\\
-7 & 4 & -5 & -1\\
-9 & 5 & -6 & 0\\
3 & -1 & 1 & \mathbf{3}
\end{array}\right|\Rightarrow\left|\begin{array}{rrr}
-15 & 2 & 1\\
-18 & \mathbf{11} & -14\\
-27 & 15 & -18
\end{array}\right|\Rightarrow\left|\begin{array}{rr}
-129 & 39\\
27 & 12
\end{array}\right|\Rightarrow-495,
\end{equation}
 hence, we have, 
\begin{equation}
A=\frac{-495}{(-1)^{3}(3)^{2}(11)^{1}}=5.
\end{equation}

\section{Mathematical Proofs}

Let us now to prove the statements of the precedent sections. As we
know, a determinant $A$, which is composed by the elements $a_{ij}$,
can be defined by the following expression: 
\begin{equation}
A=\sum_{\alpha\beta\ldots\nu=1}^{n}\varepsilon_{\alpha\beta\ldots\nu}\left(a_{1\alpha}a_{2\beta}\ldots a_{n\nu}\right),\label{definition}
\end{equation}
 where $\varepsilon_{\alpha\beta\ldots\nu}$ it is the well-known
Levi-Civita symbol -- a complete anti-symmetrical tensor with $\varepsilon_{12\ldots n}=1$.

For a future use we also present two important theorems of the theory
of determinants, the Laplace and Cauchy theorems. Let $a_{ij}$ be
an element of the determinant $A$ and let $A_{ij}$ be the minor
obtained from $A$ when we eliminate its line $i$ and column $j$.
We define thus the \textit{cofactor} of $a_{ij}$, and we represent
it by $a^{ij}$, by the expression 
\begin{equation}
a^{ij}=(-1)^{i+j}{A_{ij}}.
\end{equation}
 With the concept of cofactor we can enunciate the Laplace theorem
as follows: the summation of the elements of a given line (or column)
of the determinant $A$, multiplied each one by its respective cofactor,
equals the value of the determinant $A$. That is, 
\begin{equation}
\sum_{\alpha=1}^{n}a_{\alpha j}a^{\alpha j}=\sum_{\beta=1}^{n}a_{i\beta}a^{i\beta}=A.
\end{equation}
 At the same foot, the Cauchy theorem establishes that the summation
of a given line (or column) of the determinant $A$, each one multiplied
by the cofactors associated with the elements of any other parallel
line (or column) of $A$, it is zero. That is, 
\begin{equation}
\sum_{\alpha=1}^{n}a_{\alpha k}a^{\alpha l}=\sum_{\beta=1}^{n}a_{k\beta}a^{l\beta}=0,\hspace{1cm}k\neq l.
\end{equation}
 The Laplace and Cauchy theorem can be written in a compact form with
the definition of adjugate determinant, $A^{\dag}$. This is just
the determinant built up from the cofactors of $A$, namely, 
\begin{equation}
A^{\dag}=|a^{ij}|=|a^{ji}|.
\end{equation}
 From this definition it follows that the products $AA^{\dag}$ or
$A^{\dag}A$ result in a determinant whose diagonal elements are all
equal to $A$ (by the Laplace's theorem) and whose non-diagonal elements
vanish (by the Cauchy's theorem). From which we obtain the important
relation 
\begin{equation}
AA^{\dag}=A^{\dag}A=A^{n}.
\end{equation}

\subsection{Proof for the First Form of the Reduction Method}

Now we shall present the demonstration for the reduction method shown
on the section 3. Before this, however, it is necessary to prove the
following \emph{Lemma}: Let $A$ be a $n$-order determinant and let
$B$ be a second-order determinant whose elements $b_{ij}$ are obtained
from $A$ as below: 
\begin{eqnarray}
 & b_{11}=\left|\begin{array}{ccc}
a_{11} & \hdots & a_{1,n-1}\\
\vdots & \ddots & \vdots\\
a_{n-1,1} & \hdots & a_{n-1,n-1}
\end{array}\right|,\hspace{0.75cm}b_{12}=\left|\begin{array}{ccc}
a_{12} & \hdots & a_{1n}\\
\vdots & \ddots & \vdots\\
a_{n-1,2} & \hdots & a_{n-1,n}
\end{array}\right|,\nonumber \\
\nonumber \\
 & b_{21}=\left|\begin{array}{ccc}
a_{21} & \hdots & a_{2,n-1}\\
\vdots & \ddots & \vdots\\
a_{n1} & \hdots & a_{n,n-1}
\end{array}\right|,\hspace{0.75cm}b_{22}=\left|\begin{array}{ccc}
a_{22} & \hdots & a_{2n}\\
\vdots & \ddots & \vdots\\
a_{n2} & \hdots & a_{nn}
\end{array}\right|,
\end{eqnarray}
 In addiction, let $C$ be the central minor of $A$, whose order
is $n-2$. Namely, 
\begin{equation}
C=\left|\begin{array}{ccc}
a_{22} & \hdots & a_{2,n-1}\\
\vdots & \ddots & \vdots\\
a_{n-1,2} & \hdots & a_{n-1,n-1}
\end{array}\right|.
\end{equation}
 With these definitions, we affirm that, if $C\neq0$, then $B=AC$.
To prove this assertion, we may begin noticing that the elements of
$B$ are related to the cofactors of $A$ by the expressions 
\begin{equation}
b_{11}=a^{nn},\hspace{0.75cm}b_{12}=(-1)^{n+1}a^{n1},\hspace{0.75cm}b_{21}=(-1)^{n+1}a^{1n},\hspace{0.75cm}b_{22}=a^{11}.
\end{equation}
 Therefore, we can write 
\begin{equation}
B=\left|\begin{array}{cc}
a^{nn} & a^{n1}\\
a^{1n} & a^{11}
\end{array}\right|=\left|\begin{array}{cc}
a^{11} & a^{n1}\\
a^{1n} & a^{nn}
\end{array}\right|,
\end{equation}
 since the factor $(-1)^{n+1}$ does not matter in the calculation.
Now, consider the determinant 
\begin{equation}
D=\left|\begin{array}{ccccc}
a^{11} & 0 & \hdots & 0 & a^{n1}\\
0 & 1 & \vdots & 0 & 0\\
\vdots & \vdots & \ddots & \vdots & \vdots\\
0 & 0 & \vdots & 1 & 0\\
a^{1n} & 0 & \hdots & 0 & a^{nn}
\end{array}\right|,
\end{equation}
 which has the same value of $B$ -- as can be easily proof directly
from the definition (\ref{definition}). The determinants $D$ an
$A$ have the same order, hence we can perform the product $E=AD$.
We get, 
\begin{equation}
E=\left|\begin{array}{ccccc}
(a_{11}a^{11}+a_{1n}a^{1n}) & a_{12} & \hdots & a_{1,n-1} & (a_{11}a^{n1}+a_{1n}a^{nn})\\
(a_{21}a^{11}+a_{2n}a^{1n}) & a_{22} & \hdots & a_{2,n-1} & (a_{21}a^{n1}+a_{2n}a^{nn})\\
\vdots & \vdots & \ddots & \vdots & \vdots\\
(a_{n-1,1}a^{11}+a_{n-1,n}a^{1n}) & a_{n-1,2} & \hdots & a_{n-1,n-1} & (a_{n-1,1}a^{n1}+a_{n-1,n}a^{nn})\\
(a_{n1}a^{11}+a_{nn}a^{1n}) & a_{n2} & \hdots & a_{n,n-1} & (a_{n1}a^{n1}+a_{nn}a^{nn})
\end{array}\right|.
\end{equation}

Now, multiply the columns $2,3,\ldots,n-1$ of $E$ respectively by
$a^{12},a^{13},\ldots,a^{1,n-1}$ and add the result to its first
column. In the same way, multiply the columns $2,3,\ldots,n-1$ of
$E$ by $a^{n2},a^{n3},\ldots,a^{n,n-1}$ and add the result to its
last column. After doing this, we shall obtain a determinant $F$
whose element $f_{11}$ and $f_{nn}$ are equal to the value of the
determinant $A$ (by the Laplace's theorem) and whose others elements
of these two columns are zero (by the Cauchy's theorem). Therefore,
we have, 
\begin{equation}
F=A^{2}\left|\begin{array}{ccccc}
1 & a_{12} & \hdots & a_{1,n-1} & 0\\
0 & a_{22} & \hdots & a_{2,n-1} & 0\\
\vdots & \vdots & \ddots & \vdots & \vdots\\
0 & a_{n-1,2} & \hdots & a_{n-1,n-1} & 0\\
0 & a_{n2} & \hdots & a_{n,n-1} & 1
\end{array}\right|=A^{2}\left|\begin{array}{ccc}
a_{22} & \hdots & a_{2,n-1}\\
\vdots & \ddots & \vdots\\
a_{n-1,2} & \hdots & a_{n-1,n-1}
\end{array}\right|=A^{2}C,
\end{equation}
 where we have used the Laplace method in the last passage, expanding
$F$ by its first and last column. Therefore, we have $F=A^{2}C$,
but since $F$ was obtained from $E$ only from elementary operations
(we had added to the first and last columns of $E$ a linear combination
of the others columns), the determinant $F$ equals the determinant
$E$. Thus, we have $A^{2}C=F=E=AD=AB$, that is, $B=AC$, which proves
the \emph{Lemma}.

With the aid of this \emph{Lemma}, the first form of the reduction
method can be proven from mathematical induction. To do this, let
$A$ be a $(n+1)$-order determinant. Suppose the method is valid
for $n$-order determinants, so that each $n$-order minor of $A$,
that are built up from $n$ adjacent lines and columns of $A$, can
be reduced to a single number after $n-1$ reductions. Furthermore,
the central minor of $A$, whose order is $n-2$ (and will be called
$C$), can be reduced to a single number after $n-3$ reductions.
The reduction of that $n$-order minors will give place to four numbers
$b_{ij}$, so that, by an adequate ordination of them, we can form
the second-order determinant $B=|b_{ij}|$. So, if the reduction method
holds for the $(n+1)$-order determinant $A$ -- that is what we want
to proof --, then the determinant $A$ must be equal to the determinant
$B$ divided by the central element of the previous reduced determinant,
that is, the minor $C$. Therefore, we must have $A=B/C$ (since we
suppose that $C\neq0$). But that is precisely what we have proven
on the \emph{Lemma} above. Hence, since the reduction method is valid
for third-order determinants, the method also will hold for fourth-order
determinants, and hence for fifth-orders as well, and so on. Therefore,
the reduction method will be valid for a general $n$-order determinant.
This completes the proof.

\subsection{Proof for the Second Form of the Reduction Method}

The proof for this alternative form of the reduction method is, in
fact, much easier than the proof for the precedent one. Consider again
a $n$-order determinant $A$ and let $a_{rs}$ be a non-zero element
of it. If we fix the element $a_{rs}$ from $A$ and reduce it as
we have done in the section 4, the reduced determinant $B$ will assume
the following form: 
\begin{equation}
B=\left|\begin{array}{ccc}
(a_{11}a_{rs}-a_{1s}a_{r1}) & \hdots & -(a_{1n}a_{rs}-a_{1s}a_{rn})\\
\vdots & \ddots & \vdots\\
-(a_{n1}a_{rs}-a_{ns}a_{r1}) & \hdots & (a_{nn}a_{rs}-a_{ns}a_{rn})
\end{array}\right|.
\end{equation}

As commented on the section $4$, it is enough to prove that $A=B/(a_{rs})^{n-2}$.
To do this, we may multiply the $n-1$ lines of $A$, with $i\neq r$,
by $a_{rs}$ to obtain the determinant 
\begin{equation}
C=\left|\begin{array}{ccccc}
a_{11}a_{rs} & \ldots & a_{1s}a_{rs} & \ldots & a_{1n}a_{rs}\\
\vdots & \ddots & \vdots & \ddots & \vdots\\
a_{r1} & \ldots & a_{rs} & \ldots & a_{rn}\\
\vdots & \ddots & \vdots & \ddots & \vdots\\
a_{n1}a_{rs} & \ldots & a_{ns}a_{rs} & \ldots & a_{nn}a_{rs}
\end{array}\right|.
\end{equation}
 We have, therefore, $A=(a_{rs})^{n-1}C$. Now, we add to each element
of $C$, with $i\neq r$, the quantity $-a_{is}a_{rj}$, from which
we find 
\begin{equation}
D=\left|\begin{array}{ccccc}
a_{11}a_{rs}-a_{1s}a_{r1} & \ldots & 0 & \ldots & a_{1n}a_{rs}-a_{1s}a_{rn}\\
\vdots & \ddots & \vdots & \ddots & \vdots\\
a_{r1} & \ldots & a_{rs} & \ldots & a_{rn}\\
\vdots & \ddots & \vdots & \ddots & \vdots\\
a_{n1}a_{rs}-a_{ns}a_{r1} & \ldots & 0 & \ldots & a_{nn}a_{rs}-a_{ns}a_{rn}
\end{array}\right|,
\end{equation}
 and since we have applied only elementary operations on $D$, we
have $D=C$.

After, multiply that lines of $D$ for which $i>r$, and also that
columns for which $j>s$, by $-1$. We will find so the determinant
\begin{equation}
E=\left|\begin{array}{ccccc}
(a_{11}a_{rs}-a_{1s}a_{r1}) & \ldots & 0 & \ldots & -(a_{1n}a_{rs}-a_{1s}a_{rn})\\
\vdots & \ddots & \vdots & \ddots & \vdots\\
a_{r1} & \ldots & a_{rs} & \ldots & -a_{rn}\\
\vdots & \ddots & \vdots & \ddots & \vdots\\
-(a_{n1}a_{rs}-a_{ns}a_{r1}) & \ldots & 0 & \ldots & (a_{nn}a_{rs}-a_{ns}a_{rn})
\end{array}\right|,
\end{equation}
 and since we have multiplied $n-r$ lines, and $n-s$ columns, of
$D$ by $-1$, we get 
\begin{equation}
E=(-1)^{n-r}(-1)^{n-s}D=(-1)^{r+s}D.\label{E}
\end{equation}
 Now, if we apply the Laplace method for the column $j=s$ of the
determinant $E$, we will obtain 
\begin{equation}
E=(-1)^{r+s}a_{rs}\left|\begin{array}{ccc}
(a_{11}a_{rs}-a_{1s}a_{r1}) & \ldots & -(a_{1n}a_{rs}-a_{1s}a_{rn})\\
\vdots & \ddots & \vdots\\
-(a_{n1}a_{rs}-a_{ns}a_{r1}) & \ldots & (a_{nn}a_{rs}-a_{ns}a_{rn})
\end{array}\right|=(-1)^{r+s}a_{rs}B,\label{E2}
\end{equation}
 but since $D=C=(a_{rs})^{n-1}A$, follows from (\ref{E}) and (\ref{E2})
that $a_{rs}B=(a_{rs})^{n-1}A$, that is, 
\begin{equation}
B=(a_{rs})^{n-2}A,
\end{equation}
 which proves our proposition. Finally, since the determinant $A$
can be evaluated by a recursive use of this formula, this completes
the proof.

\section*{Addendum}

After the composition of this text I was informed that similar methods
to this presented here was already developed by two mathematicians
of the XIX century. Due to this, I wish to comment something about.

The first mathematician commented above is Mr. Charles L. Dodgson%
\footnote{C. L. Dodgson is most known by the pseudonymous Lewis Carroll -- the
famous writer of {}``\emph{Alice's Adventures in Wonderland},''
among others histories.%
}, who developed essentially the same form of the reduction method
that I had shown in the sections 2 and 3, and which he called it {}``\emph{condensation
method}'' \cite{Dodgson}. In his approach, nevertheless, it was
not presented a general proof, but only a demonstration for third
and fourth-order determinants. Moreover, the problem of the division
by zero was only attacked by the procedure of the permutation of lines
and columns of the original determinant, which is not enough in the
general, as we have seen. Nevertheless, it should also be noted that
Mr. Dodgson had used his method to develop a very interesting way
of solving linear systems of equations \cite{Dodgson}.

The other mathematician is Felice Chiò%
\footnote{F. Chiò was an Italian mathematician. He was also a disciple of Amadeo
Avogadro.%
}, who developed essentially the same form of the reduction's method
that I have shown in the section 4 \cite{Chio}. Sometimes his method
is cited as {}``\emph{pivotal method}.'' Unfortunately, I did not
have access to his original paper, so I have not commented about this.

It is of noticing, however, that the work of these two mathematicians
seems to be independent one to another and, actually, they were regarded,
until now, as independent of each other, while I showed here that
these two approaches have the same origin and can, therefore, be considered
as two different forms of a same technique.

With respect to my own research, I would like to inform that I begin
my studies on the field about ten years ago, when I was just a student
of secondary school. At that moment, I could develop only the reduction
method for third-order determinants, but, unfortunately, at that time
these questions were to me just a matter of curiosity, and I was discouraged
from them soon after. I only returned to subject when I arrived at
my graduating college, about five years ago, when I eventually browsed
through my old annotations. At this time, I generalized the reduction
method for higher-order determinants as well as I found the other
formulation of this method which was presented at the section 4. The
mathematical demonstrations took a little more of time and an appropriated
presentation of the subject was possible only recently. All work was
made in a complete independent way and I did have no knowledge of
the works of Mr. Dodgson and Mr. Chiò during this time.

Finally, I wish to point out that this technique of solving determinants
seems to be an almost unknown technique. In fact, it is not even commented
on the majority of books, specialized or not in the theory of determinants.
It is unknown by the majority of the teachers as well. Thus, I think
that the present work may provide such a popularization of the reduction
method, as well as be helpful to the teaching of the theory to students
and also to the development of the theory of determinants itself.
\selectlanguage{brazil}%
\begin{description}
\item \end{description}
\selectlanguage{english}%

\end{document}